\documentclass[12pt]{amsart}
\usepackage{amsmath,amssymb,amsbsy,amsfonts,latexsym,amsopn,amstext,cite,
                                               amsxtra,euscript,amscd,bm}
\usepackage{url}

\usepackage{mathrsfs}

\usepackage{color}
\usepackage[colorlinks,linkcolor=blue,anchorcolor=blue,citecolor=blue,backref=page]{hyperref}
\usepackage{color}
\usepackage{graphics,epsfig}
\usepackage{graphicx}
\usepackage{float}
\usepackage{epstopdf}
\hypersetup{breaklinks=true}

\usepackage[np]{numprint}
\npdecimalsign{\ensuremath{.}}

\usepackage{bibentry}

\usepackage[english]{babel}
\usepackage{mathtools}
\usepackage{todonotes}
\usepackage{url}
\usepackage[colorlinks,linkcolor=blue,anchorcolor=blue,citecolor=blue,backref=page]{hyperref}

\usepackage[norefs,nocites]{refcheck}

\begin{document}

\newcommand{\bs}{\boldsymbol}
\def \a{\alpha} \def \b{\beta} \def \d{\delta} \def \e{\varepsilon} \def \g{\gamma} \def \k{\kappa} \def \l{\lambda} \def \s{\sigma} \def \t{\theta} \def \z{\zeta}

\newcommand{\mb}{\mathbb}

\newtheorem{theorem}{Theorem}
\newtheorem{lemma}[theorem]{Lemma}
\newtheorem{claim}[theorem]{Claim}
\newtheorem{cor}[theorem]{Corollary}
\newtheorem{conj}[theorem]{Conjecture}
\newtheorem{prop}[theorem]{Proposition}
\newtheorem{definition}[theorem]{Definition}
\newtheorem{question}[theorem]{Question}
\newtheorem{example}[theorem]{Example}
\newcommand{\hh}{{{\mathrm h}}}
\newtheorem{remark}[theorem]{Remark}

\numberwithin{equation}{section}
\numberwithin{theorem}{section}
\numberwithin{table}{section}
\numberwithin{figure}{section}

\def\sssum{\mathop{\sum\!\sum\!\sum}}
\def\ssum{\mathop{\sum\ldots \sum}}
\def\iint{\mathop{\int\ldots \int}}

\newcommand{\diam}{\operatorname{diam}}

\def\squareforqed{\hbox{\rlap{$\sqcap$}$\sqcup$}}
\def\qed{\ifmmode\squareforqed\else{\unskip\nobreak\hfil
\penalty50\hskip1em \nobreak\hfil\squareforqed
\parfillskip=0pt\finalhyphendemerits=0\endgraf}\fi}

\newfont{\teneufm}{eufm10}
\newfont{\seveneufm}{eufm7}
\newfont{\fiveeufm}{eufm5}
%
%
\newfam\eufmfam
     \textfont\eufmfam=\teneufm
\scriptfont\eufmfam=\seveneufm
     \scriptscriptfont\eufmfam=\fiveeufm
%
%
\def\frak#1{{\fam\eufmfam\relax#1}}

\newcommand{\bflambda}{{\boldsymbol{\lambda}}}
\newcommand{\bfmu}{{\boldsymbol{\mu}}}
\newcommand{\bfxi}{{\boldsymbol{\eta}}}
\newcommand{\bfrho}{{\boldsymbol{\rho}}}

\def\eps{\varepsilon}

\def\fK{\mathfrak K}
\def\fT{\mathfrak{T}}
\def\fL{\mathfrak L}
\def\fR{\mathfrak R}
\def\fQ{\mathfrak Q}

\def\fA{{\mathfrak A}}
\def\fB{{\mathfrak B}}
\def\fC{{\mathfrak C}}
\def\fM{{\mathfrak M}}
\def\fS{{\mathfrak  S}}
\def\fU{{\mathfrak U}}

\def\sssum{\mathop{\sum\!\sum\!\sum}}
\def\ssum{\mathop{\sum\ldots \sum}}
\def\dsum{\mathop{\quad \sum \qquad \sum}}
\def\iint{\mathop{\int\ldots \int}}
 
\def\T {\mathsf {T}}
\def\Tor{\mathsf{T}_d}
\def\Tore{\widetilde{\mathrm{T}}_{d} }

\def\sM {\mathsf {M}}
\def\sL {\mathsf {L}}
\def\sK {\mathsf {K}}
\def\sP {\mathsf {P}}

\def\ss{\mathsf {s}}

\def \balpha{\bm{\alpha}}
\def \bbeta{\bm{\beta}}
\def \bgamma{\bm{\gamma}}
\def \bdelta{\bm{\delta}}
\def \beps{\bm{\varepsilon}}
\def \bzeta{\bm{\zeta}}
\def \blambda{\bm{\lambda}}
\def \bchi{\bm{\chi}}
\def \bphi{\bm{\varphi}}
\def \bpsi{\bm{\psi}}
\def \bxi{\bm{\xi}}
\def \bmu{\bm{\mu}}
\def \bnu{\bm{\nu}}
\def \btheta{\bm{\vartheta}}
\def \bomega{\bm{\omega}}

\def \bell{\bm{\ell}}

\def\eqref#1{(\ref{#1})}

\def\vec#1{\mathbf{#1}}

\newcommand{\abs}[1]{\left| #1 \right|}

\def\Zq{\mathbb{Z}_q}
\def\Zqx{\mathbb{Z}_q^*}
\def\Zd{\mathbb{Z}_d}
\def\Zdx{\mathbb{Z}_d^*}
\def\Zf{\mathbb{Z}_f}
\def\Zfx{\mathbb{Z}_f^*}
\def\Zp{\mathbb{Z}_p}
\def\Zpx{\mathbb{Z}_p^*}
\def\cM{\mathcal M}
\def\cE{\mathcal E}
\def\cH{\mathcal H}

\def\le{\leqslant}
\def\leq{\leqslant}
\def\ge{\geqslant}
\def\leq{\leqslant}

\def\sfB{\mathsf {B}}
\def\sfC{\mathsf {C}}
\def\sfS{\mathsf {S}}
\def\sfI{\mathsf {I}}
\def\sfT{\mathsf {T}}
\def\L{\mathsf {L}}
\def\FF{\mathsf {F}}

\def\sE {\mathscr{E}}
\def\sS {\mathscr{S}}

\def\cA{{\mathcal A}}
\def\cB{{\mathcal B}}
\def\cC{{\mathcal C}}
\def\cD{{\mathcal D}}
\def\cE{{\mathcal E}}
\def\cF{{\mathcal F}}
\def\cG{{\mathcal G}}
\def\cH{{\mathcal H}}
\def\cI{{\mathcal I}}
\def\cJ{{\mathcal J}}
\def\cK{{\mathcal K}}
\def\cL{{\mathcal L}}
\def\cM{{\mathcal M}}
\def\cN{{\mathcal N}}
\def\cO{{\mathcal O}}
\def\cP{{\mathcal P}}
\def\cQ{{\mathcal Q}}
\def\cR{{\mathcal R}}
\def\cS{{\mathcal S}}
\def\cT{{\mathcal T}}
\def\cU{{\mathcal U}}
\def\cV{{\mathcal V}}
\def\cW{{\mathcal W}}
\def\cX{{\mathcal X}}
\def\cY{{\mathcal Y}}
\def\cZ{{\mathcal Z}}
\newcommand{\rmod}[1]{\: \mbox{mod} \: #1}

\def\cg{{\mathcal g}}

\def\vy{\mathbf y}
\def\vr{\mathbf r}
\def\vx{\mathbf x}
\def\va{\mathbf a}
\def\vb{\mathbf b}
\def\vc{\mathbf c}
\def\ve{\mathbf e}
\def\vf{\mathbf f}
\def\vg{\mathbf g}
\def\vh{\mathbf h}
\def\vk{\mathbf k}
\def\vm{\mathbf m}
\def\vz{\mathbf z}
\def\vu{\mathbf u}
\def\vv{\mathbf v}

\def\e{{\mathbf{\,e}}}
\def\ep{{\mathbf{\,e}}_p}
\def\eq{{\mathbf{\,e}}_q}
\def\er{{\mathbf{\,e}}_r}
\def\es{{\mathbf{\,e}}_s}

 \def\SS{{\mathbf{S}}}

 \def\0{{\mathbf{0}}}
 
 \newcommand{\GL}{\operatorname{GL}}
\newcommand{\SL}{\operatorname{SL}}
\newcommand{\lcm}{\operatorname{lcm}}
\newcommand{\ord}{\operatorname{ord}}
\newcommand{\Tr}{\operatorname{Tr}}
\newcommand{\Span}{\operatorname{Span}}

\def\({\left(}
\def\){\right)}
\def\l|{\left|}
\def\r|{\right|}
\def\fl#1{\left\lfloor#1\right\rfloor}
\def\rf#1{\left\lceil#1\right\rceil}
\def\sumstar#1{\mathop{\sum\vphantom|^{\!\!*}\,}_{#1}}

\def\mand{\qquad \mbox{and} \qquad}

\def\tblue#1{\begin{color}{blue}{{#1}}\end{color}}




\hyphenation{re-pub-lished}

\mathsurround=1pt

\def\bfdefault{b}

\def \F{{\mathbb F}}
\def \K{{\mathbb K}}
\def \N{{\mathbb N}}
\def \Z{{\mathbb Z}}
\def \P{{\mathbb P}}
\def \Q{{\mathbb Q}}
\def \R{{\mathbb R}}
\def \C{{\mathbb C}}
\def\Fp{\F_p}
\def \fp{\Fp^*}

 \def \xbar{\overline x}

\title{On the number of Diophantine $m$-tuples in finite fields}

\author[I. E. Shparlinski] {Igor E. Shparlinski}
\address{School of Mathematics and Statistics, University of New South Wales, Sydney NSW 2052, Australia}
\email{igor.shparlinski@unsw.edu.au}

\begin{abstract}   We use a new argument to improve the error term in the asymptotic 
formula for the number of 
Diophantine $m$-tuples in finite fields, 
which is due to A.~Dujella and M.~Kazalicki  (2021) 
and  N.~Mani and  S.~Rubinstein-Salzedo  (2021). 
\end{abstract}


\maketitle

%


\section{Introduction}

\subsection{Motivation and set-up}  We recall the classical definition of a
 {\it Diophantine $m$-tuple\/} as a vector $\(a_1, \ldots, a_m\) \in \N^m$ such 
 that all shifted products $a_ia_j+1$, $1 \le i < j \le m$, are perfect squares.
 
The long-standing conjecture on the finiteness of  the set 
of  Diophantine quintuples, after a series of intermediate results by various 
authors,  has been established in a striking work of Dujella~\cite{Duj}, who has also shown the 
the non-existence of Diophantine sextuples. More recently, He, Togb{\'e} and 
Ziegler~\cite{HTZ} have show 
 the non-existence of Diophantine quintuples is shown, see also~\cite{BCM}.
  Quite naturally,  these results
have suggested to study the generalisation of this notion to other algebraic domains 
such as, for example,  the set of rational numbers or points on curves, as well as in 
many other directions, see,
for example,~\cite{DKRM, DKMS, Gup, KN} and references therein. The notion also readily extends to
the setting of finite fields, see~\cite{DK,HKMNOSS,MR-S}.

Let $q$ be an odd  prime power and let $\F_q$ be the finite field of $q$ elements. 

For $r \in \F_q^*$, we say that an $m$-tuple $\(a_1, \ldots, a_m\) \in \F_q^m$  
form a  {\it Diophantine $m$-tuple in $\F_q$  with a shift $r$} 
if all $m(m-1)/2$ shifted products $a_ia_j + r$ are perfect squares in $\F_q$. 

\begin{remark}
We note that it is customary to exclude zero values from the domain from which 
$a_1, \ldots, a_m$ are drawn. However in the counting results below this makes 
no difference, while this simplifies the notation. In particular, 
the total number  of such  $m$-tuples over $\F_q$ with a zero entry (which is at most $mq^{m-1}$)
can be absorbed in the error term of our  asymptotic formula. 
\end{remark}

Let $N_r(m,q)$ be the number of distinct Diophantine $m$-tuple in $\F_q$  with a shift $r$. 
It has been shown by Dujella and  Kazalicki~\cite{DK} that for $r=1$ and a prime $p$
we have 
 \begin{equation}
\label{eq:D-K Bound}
N_1(m,p) = 2^{-m(m-1)/2} p^m + o(p^m).
\end{equation}
Using some ideas of Dujella and  Kazalicki~\cite{DK}, 
Mani and Rubinstein-Salzedo~\cite[Theorem~5.1]{MR-S} have given explicit formulas 
for  $N_r(2,q) $ and $N_r(3,q)$ and for $m \ge 4$ presented 
a more precise than~\eqref{eq:D-K Bound} asymptotic formula
 \begin{equation}
\label{eq:MR-S Bound}
N_r(m,q) = 2^{-m(m-1)/2} q^m + O\(q^{m-1/2}\), 
\end{equation}
where the implied constant may depend on $m$, 
which also holds for any $r \in\F_q^*$ (and it is also easy to see
that for any odd prime power $q$ rather than just for a prime $q = p$ as in~\cite{MR-S}). 
We also observe that the bound~\eqref{eq:MR-S Bound} can be derived within 
the initial approach of  Dujella and  Kazalicki~\cite{DK} if one appeals to 
a version of the Lang-Weil bound~\cite{LaWe}.  

\subsection{New bound}

Here we show that using some simple arguments the bound on the error term
in~\eqref{eq:MR-S Bound}  can be improved. 

\begin{theorem} 
\label{thm:Nmq}
 For a fixed $m\ge 4$, uniformly over $r \in\F_q^*$,  we have 
$$
N_(m,q) = 2^{-m(m-1)/2} q^m + O\(q^{m-1}\), 
$$
where the implied constant may depend on $m$. 
\end{theorem}

As in~\cite{MR-S} our proof is based on an application of the {\it Weil bound\/} for multiplicative 
character sums with polynomials, see, for example,~\cite[Theorem~11.23]{IwKow}.

\section{Proof of Theorem~\ref{thm:Nmq}}  

\subsection{Preliminary transformations} 
Since there are $O\(q^{m-1}\)$ choices of $m$-tuples $\(a_1, \ldots, a_m\) \subseteq \F_q^m$  for which $a_ia_j + r = 0$ for some 
$1 \le i < j \le m$, or with $a_i = 0$  for some 
$1 \le i  \le m$, following the argument of~\cite{MR-S}, 
we write 
$$
N_r(m,q) = 2^{-m(m-1)/2} \sum_{a_1, \ldots, a_m\in  \F_q^*}
 \prod_{1 \le i < j \le m}\(1 + \chi\(a_i a_j + r\)\) + O\(q^{m-1}\), 
 $$
 where $\chi$ is the quadratic character of $\F_q$, we refer to~\cite[Chapter~3]{IwKow} for a 
 background on characters. 
 Therefore
 \begin{equation}
\label{eq:N and R}
N_r(m,q) = 2^{-m(m-1)/2} q^m + \sum_{\substack{\beps \in \{0,1\}^{m(m-1)/2}\\ \beps \ne \mathbf 0}}  
R\(\beps\)
+ O\(q^{m-1}\), 
\end{equation}
 where for $\beps = \(\varepsilon_{i,j}\)_{1 \le i < j \le m} \in \{0,1\}^m$
 \begin{equation}
\label{eq:Reps}
R\(\beps\) = \sum_{a_1, \ldots, a_m\in  \F_q^*}
 \prod_{1 \le i < j \le m}\chi\(a_i a_j + r\)^{\varepsilon_{i,j}} . 
\end{equation}
 We now fix $\beps \in \{0,1\}^{m(m-1)/2}$ with  $\beps \ne \mathbf 0$ and 
 estimate $R\(\beps\)$. Renumbering the variables $a_1, \ldots, a_m$, 
 we see that  without loss of generality, we can assume that 
 \begin{equation}
\label{eq:eps12}
\varepsilon_{1,2}=1.
\end{equation}

 We now consider the following two cases depending on vanishing and 
 non-vanishing of the exponents $\varepsilon_{i,j}$ with $2 \le i < j \le m$. 
 
 \subsection{Vanishing exponents $\varepsilon_{i,j}$ with $2 \le i < j \le m$} Assume that 
    \begin{equation}
\label{eq:vanish}
 \varepsilon_{i,j} = 0, \qquad  \text{for all $2 \le i < j \le m$.}
 \end{equation}Then we see that under the conditions~\eqref{eq:vanish}
 the expression for $R\(\beps\)$ in~\eqref{eq:Reps} simplifies as 
\begin{align*}
R\(\beps\) & = \sum_{a_1, \ldots, a_m\in  \F_q^*}
 \prod_{2 \le  j \le m}\chi\(a_1 a_j + r\)^{\varepsilon_{1,j}} \\
 & = \sum_{a_1 \in  \F_q^*}  \prod_{2 \le  j \le m}  \sum_{a_j \in  \F_q^*} \chi\(a_1 a_j + r\)^{\varepsilon_{1,j}}.
\end{align*}
Hence, estimating the sums over $a_3, \ldots, a_m$ trivially as $q-1$ and recalling 
our assumption~\eqref{eq:eps12},  we obtain 
$$
|R\(\beps\)|   \le (q-1)^{m-2} 
\sum_{a_1 \in  \F_q^*}   \left| \sum_{a_2 \in  \F_q^*} \chi\(a_1 a_2 + r\)\right| . 
$$
Clearly, for every $a_1 \in  \F_q^*$ we have 
\begin{align*}
 \sum_{a_2 \in  \F_q^*} \chi\(a_1 a_2 + r\) & =  \sum_{a \in  \F_q^*} \chi\(a + r\)\\
 & = \sum_{a \in  \F_q} \chi\(a \) - \chi(1) =   \sum_{a \in  \F_q} \chi\(a\) - 1=-1.
\end{align*}
 Hence we obtain 
  \begin{equation}
\label{eq:Bound I}
|R\(\beps\)|  \le  (q-1)^{m-1} 
\end{equation}
in this case. 

 \subsection{Non-vanishing exponents $\varepsilon_{i,j}$ with $2 \le i < j \le m$} We now assume that 
   \begin{equation}
\label{eq:non-vanish}
 \varepsilon_{i,j} \ne 0, \qquad  \text{for some $2 \le i < j \le m$.}
 \end{equation}
 We write $R\(\beps\)$ as 
 $$
 R\(\beps\)
 =  \sum_{a_1, \ldots, a_m\in  \F_q^*}  \prod_{2 \le  j \le m}\chi\(a_1 a_j + r\)^{\varepsilon_{1,j}}  \prod_{2\le i < j \le m} \chi\(a_i a_j + r\)^{\varepsilon_{i,j}} . 
 $$ 
 
 Observe that for any $b \in \F_q^*$ the map
$$(a_1, a_2, … , a_m) \mapsto (a_1/b, a_2b, … , a_mb)$$ 
is a permutation on $\F_p^m$. Hence 
\begin{align*}
 R\(\beps\)
  & = (p-1)^{-1} \sum_{b \in  \F_q^*}  \sum_{a_1, \ldots, a_m\in  \F_q^*}
 \prod_{2 \le  j \le m}\chi\(a_1 a_j + r\)^{\varepsilon_{1,j}} \\
 &\qquad \qquad \qquad \qquad \qquad \qquad 
  \prod_{2\le i < j \le m} \chi\(a_i a_j b^2 + r\)^{\varepsilon_{i,j}}, 
\end{align*}
which we now rearrange as 
$$
 R\(\beps\)
  = (q-1)^{-1}   \sum_{a_2, \ldots, a_m\in  \F_q^*}
  S(a_2,…, a_m)   T(a_2,…, a_m), 
  $$
where 
\begin{align*}
& S(a_2,…, a_m) = \sum_{a_1 \in  \F_q^*}  
\chi\(\prod_{2 \le  j \le m}\(a_1 a_j + r\)^{\varepsilon_{1,j}}\),\\
& T(a_2,…, a_m) =  \sum_{b\in  \F_q^*} 
\chi\( \prod_{2\le i < j \le m}\(a_i a_j b^2 + r\)^{\varepsilon_{i,j}}\).
\end{align*}

We now examine the polynomials
\begin{align*}
& F_{a_2,…, a_m}(X) = \prod_{2 \le  j \le m}\(a_j X + r\)^{\varepsilon_{1,j}}, \\
& G_{a_2,…, a_m}(X) = \prod_{2\le i < j \le m}\(a_i a_j X^2 + r\)^{\varepsilon_{i,j}}.
\end{align*}
Because of our assumptions~\eqref{eq:eps12} and~\eqref{eq:non-vanish}
both these polynomials are of positive degree.

Furthermore, it is clear that there are at most $O\(q^{m-2}\)$ choices for $(m-1)$-tuples $(a_2,…, a_m)\in \F_q^{m-1}$ 
for which at least one of the polynomials $ F_{a_2,…, a_m}(X)$ and  $G_{a_2,…, a_m}(X)$
is a perfect square in the algebraic closure of $\F_q$. In this case we estimate 
both sums $S(a_2,…, a_m) $ and $T(a_2,…, a_m)$ trivially as $q-1$. 
Hence,  the contribution to  $R\(\beps\)$ from such sums is 
   \begin{equation}
\label{eq:bad a_i}
A =  O\((q-1)^{-1}  q^{m-2} (q-1)^2\) = O(q^{m-1}).
\end{equation}
For other choices  of $(a_2,…, a_m)\in \F_q^{m-1}$, by the Weil bound, see, for example,~\cite[Theorem~11.23]{IwKow}, we have 
$$
S(a_2,…, a_m), T(a_2,…, a_m) = O(q^{1/2}).
$$
Hence,  the contribution to  $R\(\beps\)$ from such sums is 
   \begin{equation}
\label{eq:good a_i}
B =  O\((q-1)^{-1}  q^{m-1} \(q^{1/2}\)^2\) = O(q^{m-1}).
\end{equation}
Combining~\eqref{eq:bad a_i} and~\eqref{eq:good a_i}
we arrive to 
   \begin{equation}
\label{eq:Bound II}
R\(\beps\)   = A + B = O(q^{m-1}) 
\end{equation}
in this case. 

\subsection{Concluding the proof} Substituting the bounds~\eqref{eq:Bound I} 
and~\eqref{eq:Bound II} in~\eqref{eq:N and R} we immediately obtain the desired result. 

\section{Comments}

It is easy to see that all implied constants can be evaluated explicitly. 
Hence one can use our argument to estimate the smallest $q$ (in terms of $m$) 
for which $N_r(m,q) > 0$ for all $r \in \F_q^*$. However the inductive approach of 
 Dujella and Kazalicki~\cite[Theorem~17]{DK} seems to be more effective for this question.
 
 Since we have multivariate character sums, it is also natural to try improve Theorem~\ref{thm:Nmq}
 via  the use of some version of  the {\it Deligne bound\/}, see, for example~\cite{Katz}. 
 Unfortunately, our polynomials have a high dimensional 
 singularity locus, which seems to prevent this approach.
 
 \section*{Acknowledgement}  
The authors would like to thank Andrej Dujella for several very useful comments and suggestions.

This work was supported by  he Australian Research Council (Discovery Project DP200100355).

 \end{document}